\documentclass[a4paper,12pt]{amsart}
\usepackage{amssymb}
\usepackage{ifthen}
\usepackage{graphicx}
\usepackage{float}
\usepackage{caption}
\usepackage{subcaption}
\usepackage{cite}
\usepackage{amsfonts}
\usepackage{amscd}
\usepackage{amsxtra}
\usepackage{color}

\newtheorem{thm}{Theorem}
\newtheorem{cor}{Corollary}
\newtheorem{lem}{Lemma}
\newtheorem{rem}{Remark}

\newtheorem{conj}{Conjecture}
\newtheorem{prob}{Problem}

\theoremstyle{definition}
\newtheorem{defn}{Definition}[section]
\newtheorem{example}{Example}

\newenvironment{pf}[1][]{%
 \vskip 1mm
 \noindent
 \ifthenelse{\equal{#1}{}}%
  {{\slshape Proof. }}%
  {{\slshape #1.} }%
 }%
{\qed\bigskip}

\newcounter{alphabet}

\newenvironment{Thm}[1][]{\refstepcounter{alphabet}%
\bigskip%
\noindent%
{\bf Theorem \Alph{alphabet}}%
\ifthenelse{\equal{#1}{}}{}{ (#1)}%
{\bf .} \itshape}{\vskip 8pt}

\newcommand{\ID}{{\mathbb D}}

\newcommand{\dist}{{\operatorname{dist}}}

\def\be{\begin{equation}}
\def\ee{\end{equation}}

\newcommand{\bee}{\begin{enumerate}}
\newcommand{\eee}{\end{enumerate}}

\newcommand{\blem}{\begin{lem}}
\newcommand{\elem}{\end{lem}}
\newcommand{\bthm}{\begin{thm}}
\newcommand{\ethm}{\end{thm}}
\newcommand{\bcor}{\begin{cor}}
\newcommand{\ecor}{\end{cor}}
\newcommand{\beg}{\begin{example}}
\newcommand{\eeg}{\end{example}}
\newcommand{\begs}{\begin{examples}}
\newcommand{\eegs}{\end{examples}}
\newcommand{\bdefe}{\begin{defn}}
\newcommand{\edefe}{\end{defn}}
\newcommand{\bprob}{\begin{prob}}
\newcommand{\eprob}{\end{prob}}
\newcommand{\bques}{\begin{ques}}
\newcommand{\eques}{\end{ques}}
\newcommand{\bei}{\begin{itemize}}
\newcommand{\eei}{\end{itemize}}
\newcommand{\bcon}{\begin{conj}}
\newcommand{\econ}{\end{conj}}
\newcommand{\bcons}{\begin{conjs}}
\newcommand{\econs}{\end{conjs}}
\newcommand{\bprop}{\begin{propo}}
\newcommand{\eprop}{\end{propo}}
\newcommand{\br}{\begin{rem}}
\newcommand{\er}{\end{rem}}
\newcommand{\brs}{\begin{rems}}
\newcommand{\ers}{\end{rems}}
\newcommand{\bo}{\begin{obser}}
\newcommand{\eo}{\end{obser}}
\newcommand{\bos}{\begin{obsers}}
\newcommand{\eos}{\end{obsers}}
\newcommand{\bpf}{\begin{pf}}
\newcommand{\epf}{\end{pf}}
\newcommand{\ba}{\begin{array}}
\newcommand{\ea}{\end{array}}
\newcommand{\beq}{\begin{eqnarray}}
\newcommand{\beqq}{\begin{eqnarray*}}
\newcommand{\eeq}{\end{eqnarray}}
\newcommand{\eeqq}{\end{eqnarray*}}

\newcommand{\ds}{\displaystyle}

\newcounter{minutes}
\divide\time by 60
\newcounter{hours}
\multiply\time by 60 \addtocounter{minutes}{-\time}
\begin{document}

\bibliographystyle{amsplain}

%

\title[Improved Bohr's phenomenon in quasi-subordination classes]{Improved Bohr's phenomenon in quasi-subordination classes}

\thanks{
File:~\jobname .tex,
          printed: \number\day-\number\month-\number\year,
          \thehours.\ifnum\theminutes<10{0}\fi\theminutes}

\author[S. Ponnusamy]{Saminathan Ponnusamy
}
\address{
S. Ponnusamy, Department of Mathematics,
Indian Institute of Technology Madras, Chennai-600 036, India.
}
\email{samy@iitm.ac.in}

\author[R. Vijayakumar]{Ramakrishnan Vijayakumar}
\address{
R. Vijayakumar, Department of Mathematics,
Indian Institute of Technology Madras, Chennai-600 036, India.
}
\email{mathesvijay8@gmail.com}

\author[K.-J. Wirths]{Karl-Joachim Wirths}
\address{K.-J. Wirths, Institut f\"ur Analysis und Algebra, TU Braunschweig,
38106 Braunschweig, Germany.}
\email{kjwirths@tu-bs.de}


\begin{abstract}
Recently the present authors established refined versions of Bohr's inequality in the case of bounded analytic functions.
In this article, we state and prove a generalization of these results
in a reformulated ``distance form" version and thereby we extend the refined versions of the Bohr inequality for the class of the quasi-subordinations which contains both the classes of majorization and subordination as special cases. As a consequence, we obtain
several new results.
\end{abstract}

\maketitle
\pagestyle{myheadings}
\markboth{S. Ponnusamy, R. Vijayakumar and K.-J. Wirths}{Bohr's phenomenon in quasi-subordination classes}

\section{Introduction and two Main results}

Let $\mathbb{D}=\{z\in \mathbb{C}:\,|z|<1\}$ denote the open unit disk, and $\overline{\ID} =\ID\cup \partial \ID
=\{z:\,|z|\leq1\}.$
Then the classical Bohr inequality \cite{Bohr-14}, compiled by Hardy in 1914 from his correspondence with Bohr,
states the following.

\begin{Thm}\label{PVW1-theA}
If $f(z)=\sum_{n=0}^{\infty} a_n z^n$ is analytic in $\mathbb{D}$ with values in $\overline{\ID}$,
then
\be\label{KKP1-eq1}
M_f(r):=\sum_{n=0}^{\infty} |a_n|\, r^n \leq 1  
~\mbox{ for each  $r\leq 1/3$}
\ee
and  the constant $1/3$ cannot be improved.
\end{Thm}

Bohr originally proved the inequality \eqref{KKP1-eq1} only for $r\leq 1/6$ and the value $1/3$
was obtained independently by M. Riesz, I. Schur and N. Wiener. Some other proofs of this inequality \eqref{KKP1-eq1} were given
by Sidon \cite{Sidon-27-15} and Tomi\'c \cite{Tomic-62-16}. Several extensions of Theorem~A 
may be obtained from \cite{BoasKhavin-97-4,Bohr-14, Bom-62,BomBor-04,DjaRaman-2000}.
For a detailed account of literature on this topic, we refer to Abu-Muhanna et al. \cite{AAPon1}, Defant and Prengel \cite{DePre06},
Garcia et al. \cite{GarMasRoss-2018}. See also  recent works from
\cite{AliBarSoly, AlkKayPon1, BhowDas-98, DeGarM04, EvPoRa-2017, KayPon1, KayPon3, KayPon2, KayPonShak1, LiuShangXu-89}
and the references therein. Surprisingly, in a recent paper, the present authors in \cite{PVW1-preprint} refined the Bohr inequality in the following improved form.

\begin{Thm}\label{PVW1-theA-a}
Suppose that $f(z) = \sum_{n=0}^\infty a_n z^n$ is an analytic function in $\ID$, $|f(z)| \leq 1$ in $\ID$,
$f_0(z)=f(z)-a_0$, and  $\|f_0\|_r$ denotes the quantity defined by
$$\|f_0\|_r=\sum_{n=1}^{\infty}|a_n|^2r^{2n}.
$$
Then 
\begin{equation}\label{PVW1-eq1}
|a_0|+ \sum_{n=1}^\infty |a_n|r^n+  \left(\frac{1}{1+|a_0|}+\frac{r}{1-r}\right)\|f_0\|_r \leq 1 , ~\mbox{ for every }~ r \leq \frac{1}{2+|a_0|}
\end{equation}
and the numbers $\frac{1}{2+|a_0|}$ and $\frac{1}{1+|a_0|}$ cannot be improved.  Moreover,
\begin{equation}\label{PVW1-eq2}
|a_0|^2+ \sum_{n=1}^\infty |a_n|r^n+ \left(\frac{1}{1+|a_0|}+\frac{r}{1-r}\right)\|f_0\|_r \leq 1 ~\mbox{ for every  }~ r \leq \frac{1}{2}
\end{equation}
and the numbers $\frac{1}{2}$ and $\frac{1}{1+|a_0|}$ cannot be improved.
\end{Thm}

It is important to point out that $\frac{1}{3}\leq \frac{1}{2+|a_0|}\leq \frac{1}{2}$ and $1/3$ is achieved when $|a_0|=1$.
In the case $a_0=0$, we have a sharp result in \cite[Theorem 2]{PVW1-preprint}.

\br
If the constant term $|a_0|$ in \eqref{PVW1-eq1} is replaced by $|a_0|^p$ for {\color{red} $0<p\leq 2$,} then it can be easily seen from the hypothesis of
Theorem~B 
that the sharp inequality
\begin{equation}\label{PVW1-eq2e}
|a_0|^p+ \sum_{n=1}^\infty |a_n|r^n+ \left(\frac{1}{1+|a_0|}+\frac{r}{1-r}\right)\|f_0\|_r \leq 1 ~\mbox{ for every  }~ r \leq \frac{1-|a_0|^p}{2-|a_0|^2-|a_0|^p}
\end{equation}
holds, where
$$ \inf _{|a_0|<1} \left \{\frac{1-|a_0|^p}{2-|a_0|^2-|a_0|^p} \right \}=\frac{p}{2+p}.
$$
The cases $p=1,2$ obviously lead to \eqref{PVW1-eq1} and \eqref{PVW1-eq2}, respectively. The inequality \eqref{PVW1-eq2e}
follows from the proof of Theorem~B 
in \cite{PVW1-preprint}. Indeed in the proof of \cite[Theorem 1]{PVW1-preprint}, we just need to consider
$$  \Psi _p (r) =|a_0|^p+ \frac{r}{1-r} \left(1-|a_0|^2\right)
$$
and observe that $\Psi _p (r)\leq 1$ if and only if
$$r\leq \varphi (x)= \frac{1-x^p}{2-x^2-x^p}, \quad x=|a_0|\in [0,1).
$$
Moreover, for {\color{red} $0<p\leq 2$,} it is a simple exercise to see that
$$\varphi '(x)= \frac{x^{p-1}A(x)}{(2-x^2-x^p)^2}, \quad A(x) =-p-(2-p)x^2 +2x^{2-p}.
$$
Because $A'(x) =2(2-p)x^{1-p}(1-x^{p}) \geq 0 $ for {\color{red} $0< p\leq 2$} and $x\in [0,1]$, it follows that $A(x)\leq  A(1)=0 $ and thus, $\varphi$ is \
{\color{red} decreasing} on $[0,1)$. This gives
$$\varphi (x) \geq \lim_{x\rightarrow 1^-}\varphi (x) =\frac{p}{2+p}.
$$
For the sharpness of the radius in question in \eqref{PVW1-eq2e}, we consider the function $f=\varphi_a$ given by
$$\varphi_a (z)  = \frac{a-z}{1-az} =a - (1-a^2)\sum_{k=1}^\infty a^{k-1} z^{k}, \quad z\in\ID,
$$
where $a\in (0,1)$. For this function as in \cite{PVW1-preprint}, it follows that
$$\ds -a +a^p + M_{\varphi_a }(r)+ \left(\frac{1}{1+a}+\frac{r}{1-r}\right)\|\varphi_a -a\|_r  =1-a +a^p+\frac{(1- a)[(2+a)r-1]}{1-r}
$$
%
which is bigger than $1$ if and only if
$$r\leq \varphi (a)= \frac{1-a^p}{2-a^2-a^p}, \quad a\in [0,1),
 $$
and allowing $a \rightarrow 1^-$ one also gets the value $p/(2+p)$, independent of $a$.
\er

Our main concern in this article is to deal with few other related questions about the Bohr inequality. For example,  it is well-known that the
Bohr radius $1/3$ continues to hold in Theorem~A 
even if the assumption on $f$
is replaced by the condition ${\rm Re\,}f(z)<1$ in $\ID$ and $a_0=f(0)\in  [0,1)$. In fact, this condition implies that (see \cite[Carath\'{e}odory's Lemma, p.41]{DurenUniv-83-8})
$|a_n|\leq 2(1-a_0)$ for all $n\geq 1$
and thus,  we have the following sharp inequality as observed in \cite{PaulPopeSingh-02-10}
$$ M_f(r)= a_0+\sum_{n=1}^{\infty} |a_n|\, r^n \leq a_0 +2(1-a_0)\frac{r}{1-r}\leq 1 ~\mbox{  for each $r\leq 1/3$}.
$$
Therefore a natural question is to look for the analog of the refined version of it in the settings of Theorem~B. 
We answer this question in the following statement whose proof will be given in Section \ref{Sec3-PVW2}. 

\bthm\label{PVW1-cor1c}
Let $f(z)$ be an analytic function in $\ID$ such that $f(z)=\sum_{n=0}^\infty a_nz^n$, $a_{0}\in (0,1)$, and ${\rm Re}\, f(z)<1$ in  $\ID$.
Then
$$\sum_{n=0}^\infty |a_n|r^n+\left(\frac{1}{1+a_0}+\frac{r}{1-r}\right)\sum_{n=1}^{\infty}|a_n|^2r^{2n}  \le 1
$$
holds for all  $r\leq r_{*}$, where $r_{*}\approx 0.24683$ is the unique root of the equation $3r^3-5r^2-3r+1=0$ in the interval $(0,1)$.
Moreover, for any  $a_{0}\in (0,1)$ there exists a uniquely defined $r_0=r_0(a_0)\in \left(r_{*},\frac{1}{3}\right)$ such that
$$\sum_{n=0}^\infty |a_n|r^n+\left(\frac{1}{1+a_0}+\frac{r}{1-r}\right)\sum_{n=1}^{\infty}|a_n|^2r^{2n}  \le 1
$$
for $r\in [0,r_0]$. The radius $r_0=r_0(a_0)$ can be calculated as the solution of the equation
\be\label{PVW1-eq2-e}
 \Phi(\lambda,r)= 4r^3\lambda^2- (7r^3+3r^2-3r+1)\lambda +6r^3-2r^2-6r+2 = 0,
\ee
where $\lambda = 1 - a_0.$ The result is sharp.
\ethm

In Section \ref{Sec3-PVW2}, we generalize Theorem~B 
for a general class of quasi-subordinations which contains both subordination and  majorization. Furthermore, we present few other important consequences including the proof of
Theorem \ref{PVW1-cor1c}.  In Section \ref{Sec4-PVW}, we introduce Bohr's phenomenon in a refined formulation in a more general
family of subordinations.

\section{Quasi-subordination and the proof of Theorem \ref{PVW1-cor1c}} \label{Sec3-PVW2}

For any two analytic functions $f$ and $g$ in $\ID$, we say that the function $f$ is \textit{quasi-subordinate} to $g$ (relative to $\Phi$), denoted by $f(z)\prec _q g(z)$ (relative to $\Phi$) in $\ID$, if there exist two functions $\Phi $  and  $\omega$, analytic in $\ID$, satisfying $\omega(0)=0$, $|\Phi(z)|\leq 1$ and $|\omega (z)|\leq 1$ for $|z|<1$ such that
\begin{equation}\label{Eq7a}
f(z)=\Phi(z)g(\omega(z)).
\end{equation}

The case $f$ is \textit{quasi-subordinate} to $g$ (relative to $\Phi \equiv 1$)  corresponds to subordination. That is $f(z)\prec _q g(z)$ (relative to $\Phi \equiv 1$)
in $\ID$ is equivalent to saying that $f(z)\prec g(z)$, the usual subordination.
Similarly, the case $\omega(z) = z$ gives majorization, i.e. \eqref{Eq7a} reduces to the form $f(z)=\Phi (z)g(z)$.
Thus, the notion of quasi-subordination includes both the concept of subordination and the principle of majorization.
See \cite{MacGre-67, Robert-70, Rogo-43} and  the recent paper \cite{AlkKayPon1} in connection with Bohr's radius.

\subsection{Bohr's phenomenon for the class of quasi-subordinations}
For the proof of Theorem \ref{PVW1-cor1c}, we need some preparation.

\blem\label{PVW2-quasiTh}
Let $f(z)$ and $g(z)$ be two analytic functions in $\ID$ with the Taylor series expansions
$f(z)=\sum_{n = 0}^\infty a_nz^n$ and $g(z)=\sum_{n = 0}^\infty b_nz^n$ for $z\in\ID$. Suppose that
$f_0(z)=f(z)-a_0$, $g_0(z)=g(z)-a_0$ and $\|f_0\|_r$ is defined as in Theorem~B. 
If $f(z)\prec_q g(z)$ (relative to $\Phi$) then

\vspace{8pt}

$\ds \sum_{n=0}^\infty |a_n|r^n+\left(\frac{1}{1+|a_0|}+\frac{r}{1-r}\right)\|f_0\|_r   \le $
$$\sum_{n=0}^\infty|b_n|r^n +\left(\frac{1}{1+|b_0\Phi_0|}+\frac{r}{1-r}\right)(|b_0|^2(1- |\Phi_0|^2)+\|g_0\|_r)
$$
holds for all  $r\leq 1/3$, 
where $a_0=\Phi_0 b_0$ with $\Phi_0=\Phi (0)$.
\elem
\bpf 
We remark that this theorem was proved in \cite{AlkKayPon1}
without the second term on both sides of the last inequality.
Suppose that  $f\prec_q g$. Then there exist two analytic
functions $\Phi$ and $\omega$ satisfying $\omega(0)=0$, $|\omega(z)|\le1$ and
$|\Phi(z)|\le1$ for all $z\in \ID$ such that 
\begin{equation}\label{Eq7}
f(z)=\Phi(z)g(\omega(z)).
\end{equation}
Setting $z=0$ in \eqref{Eq7} gives that $a_0=\Phi_0 b_0$.
According to \cite[Theorem 2.1]{AlkKayPon1}, we obtain that
\be\label{Eq7-a}
M_f(r)= \sum_{n=0}^\infty |a_n|r^n \le  M_g(r)=\sum_{n=0}^\infty|b_n|r^n ~\mbox{ for $r\leq 1/3$.}
\ee
Finally, by \eqref{Eq7}, it follows that
$$|f(z)|^2 \leq |g(\omega(z))|^2 ~\mbox{ for }~z\in \ID
$$
and thus, as in the proof of Rogosinski's Theorem \cite{Rogo-43}, we can easily obtain that
\be\label{Eq7-b}
\|f\|_r =\sum_{n=0}^{\infty}|a_n|^2r^{2n}  \le \|g\|_r =\sum_{n=0}^{\infty}|b_n|^2r^{2n} ~\text{ for all}\ \ r \in [0,1)
\ee
and therefore, since  $a_0=\Phi_0 b_0$, we have
\be\label{Eq7-c}
\|f_0\|_r    \le |b_0|^2(1- |\Phi_0|^2)+\|g_0\|_r   ~\mbox{ for all $ r \in [0,1)$.}
\ee
The desired inequality follows from \eqref{Eq7-a}, \eqref{Eq7-b} and \eqref{Eq7-c}.
%
\epf

The following result is regarded as a generalization of Theorem~B 
and can be used to cover many situations.
Because of its independent interest, we state it here.

\blem\label{PVW1-cor1a}
Let $f(z)$ and $g(z)$ be two analytic functions in $\ID$ with the Taylor series expansions
$f(z)=\sum_{n = 0}^\infty a_nz^n$ and $g(z)=\sum_{n = 0}^\infty b_nz^n$ for $z\in\ID$. If $f(z)\prec g(z)$  then
$$\sum_{n=0}^\infty |a_n|r^n+\left(\frac{1}{1+|a_0|}+\frac{r}{1-r}\right)\sum_{n=1}^{\infty}|a_n|^2r^{2n}  \le \sum_{n=0}^\infty |b_n| r^n +\left(\frac{1}{1+|a_0|}+\frac{r}{1-r}\right)\sum_{n=1}^{\infty}|b_n|^2r^{2n}
$$  holds for all  $r\leq  1/3$.
\elem
\bpf
Set $\Phi (z)\equiv 1$. Then $\Phi_0=1$ and $a_0=b_0$. 
\epf

\bprob
Determine sharp radii in Lemmas \ref{PVW2-quasiTh} and \ref{PVW1-cor1a}.
\eprob

\subsection{Proof of Theorem \ref{PVW1-cor1c}}
Since ${\rm Re}\, f(z)<1$, we may write the given condition as
$$f(z)\prec g(z), \quad g(z) = a_{0} -2(1-a_{0})\frac{z}{1-z} =a_0 -2(1-a_{0})\sum_{n=0}^{\infty} z^{n}.
$$
Here $g(z)$ is a univalent mapping of $\ID$ onto the left half-plane  $\{ w:\, {\rm Re}\,(w)< 1\}$.
According to Lemma \ref{PVW1-cor1a}, with $g(z)=\sum_{n=0}^{\infty} b_{n} z^{n}$, it suffices to show that
$$S_{g}(r):=\sum_{n=0}^\infty |b_n| r^n +\left(\frac{1}{1+|b_0|}+\frac{r}{1-r}\right)\sum_{n=1}^{\infty}|b_n|^2r^{2n} \leq 1
~\mbox{ for every $ r \leq r_{*} $ },
$$
where $r_{*}$ is as in the statement. For convenience, we let $1-a_{0}=\lambda$ so that $a_{0}=1-\lambda$
and $b_{n}=-2\lambda$ for $n\geq 1$. This gives for $\lambda \in [0,1]$ and $r\in (0,1)$ that
\beqq
S_{g}(r)&=& 1-\lambda + 2\lambda \sum_{n=1}^\infty r^n +\left(\frac{1}{2-\lambda}+\frac{r}{1-r}\right) 4\lambda^2 \sum_{n=1}^\infty r^{2n} \\
 &=&1-\lambda \left[1 - \frac{2r}{1-r}- \left(\frac{1+r- \lambda r}{(2-\lambda)(1-r)}\right)\frac{4\lambda r^2}{1-r^2}\right]  \\
 &=&1-\lambda \left[\frac{1-3r}{1-r}-\frac{4\lambda r^2 +4 \lambda r^3 - 4\lambda^2 r^3}{(2-\lambda)(1-r)(1-r^2)}\right] \\
 &=&1-\lambda \left[\frac{1-3r}{1-r}-\frac{4\lambda r^2\{1 +(1 - \lambda)r\}}{(2-\lambda)(1-r)(1-r^2)}\right] \\
 &=&1- \lambda\left[\dfrac{\Phi(\lambda,r)}{(2-\lambda)(1-r)(1-r^2)}\right],
 \eeqq
which shows that the left hand side is less than or equal to $1$ whenever $\Phi(\lambda,r) \geq 0 ,$
where
$$ \Phi(\lambda,r)= 4r^3\lambda^2- (7r^3+3r^2-3r+1)\lambda +6r^3-2r^2-6r+2 .
$$
Before we continue, we observe from the fourth equality in the above equalities that $S_{g}(r)>1$ for $r>1/3$ and for each $\lambda \in (0,1]$.

We claim that $\Phi(\lambda,r) \geq 0$ for every $r \leq r_{*}$ and for $\lambda \in (0,1]$.
It follows that
$$ \frac{\partial^2 \Phi(\lambda,r) }{\partial \lambda^2} \geq 0 ~\mbox{ for every $\lambda \in (0,1]$}
$$
and thus,
$\frac{\partial \Phi}{\partial \lambda}$ is an increasing function of $\lambda$. This gives
$$\dfrac{\partial \Phi(\lambda,r)}{\partial \lambda} \leq  \frac{\partial \Phi}{\partial \lambda}(1,r)  = r^3-3r^2+3r-1 =-(1-r)^3,
$$
whence $\Phi$ is a decreasing function of $\lambda$ on $(0,1]$ so that
$$ \Phi(\lambda, r) \geq \Phi(1,r)=3r^3-5r^2-3r+1,
$$
which is greater than or equal to $0$ for all $r \leq r_{*}$, where $r_{*}$ is the unique root of the equation
$3r^3-5r^2-3r+1 =0$, which lies in $(0,1)$.
It is easy to see that $\Phi(1,r)=3r^3-5r^2-3r+1$ is an increasing function of $r$ in $[0,1]$ and using Mathematica
or by numerical computation by Cardano's formula, one can find that
$$r_{*}=\frac{5}{9} - \frac{2}{9} \sqrt{13} \cos\left[\frac{1}{3} ~\arctan\left(\frac{9 \sqrt{303}}{103}\right)\right]
+ \frac{2}{3}\sqrt{\frac{13}{3}} \sin\left[\frac{1}{3} ~\arctan\left (\frac{9 \sqrt{303}}{103}\right) \right]
$$
which is approximatively $0.24683$.

Since $\Phi(0,r)=2(1-3r)(1-r^2)$, we have $\Phi(0,r) \geq 0$ for $r \leq 1/3$ and  $\Phi(0,r) < 0$ for  $r < 1/3$.

Furthermore, $\Phi(1,r)=3r^3-5r^2-3r+1$ and $\Phi'(1,r) =-9r(1-r)-(r+3) < 0$ imply $\Phi(1,r)  \geq 0$ for $r \geq r_{*}$ and $\Phi(1,r) <  0$ for $r < r_{*}$. According to the fact that $ \Phi(\lambda ,r)$ is a monotonic decreasing function of $\lambda$, we see that for any $r \in (r_{*}, 1/3)$ there is a uniquely defined $\lambda(r)\in (0,1)$ such that $\Phi(\lambda(r), r) = 0.$

To prove the last assertion, we have to show that $\frac{d\lambda(r)}{dr} < 0.$ Since
$$\frac{d\lambda(r)}{dr} = -\frac{\frac{\partial \Phi(\lambda(r), r)}{\partial r}}{\frac{\partial \Phi(\lambda(r), r)}{\partial \lambda}},
$$
it is sufficient to prove that
$$\frac{\partial \Phi(\lambda(r), r)}{\partial r} <0
$$
for $\lambda \in (0,1)$ and $r \in \left[r_{*}, \frac{1}{3}\right],$ where
\beqq
\frac{\partial \Phi(\lambda(r), r)}{\partial r}
&=& 12r^2\lambda^2 - (21 r^2 + 6r -3)\lambda +18 r^2 -4r -6 \\
&=& [12r^2\lambda^2 - (21 r^2 + 6r)\lambda +3\lambda -3] + [18 r^2 -4r -3]\\
&=& -[12r^2\lambda(1-\lambda)  + 3r(3r+2)\lambda +3(1-\lambda)] - [2(1-9r^2)+4r+1].
\eeqq
This is clearly negative for the intervals in question.
%
This completes the proof of our theorem.
\hfill $\Box$


%
%
%

\section{Bohr's phenomenon for a family of subordinations} \label{Sec4-PVW}

We now turn to a discussion of Bohr's phenomenon in a refined formulation in a more general family of subordinations.
Let us first rewrite the refined version of the Bohr inequality \eqref{PVW1-eq1} in an equivalent form
$$ \sum_{n=1}^\infty |a_n|r^n+  \left(\frac{1}{2-(1-|f(0)|)}+\frac{r}{1-r}\right)\|f_0\|_r \leq   1-|a_0|=1-|f(0)|,
$$
where $\|f_0\|_r$ is defined as in Theorem~B. 
We observe that the number $1-|f(0)|$ is the distance from the point
$f(0)$ to the boundary $\partial \mathbb{D}$ of the unit disk $\mathbb{D}$ and
thus, we use this ``distance form" formulation to generalize the concept of the Bohr radius for
the class of functions $f$ analytic in $\mathbb{D}$ which take values in a given simply connected domain $\Omega$
(see also \cite{Abu}).

Now for a given univalent function $g,$ let $S(g)=\{f:\, f\prec g\}$, $\Omega =g(\ID)$
and $\dist (c,\partial \Omega)$ denote the Euclidean distance from a point $c\in \Omega$ to
the boundary $\partial \Omega$.
We say that the family $S(g)$ has a Bohr phenomenon in the refined formulation if there exists an $r_{g}$, $0<r_{g}\leq 1$,
such that whenever $f(z)=\sum_{n=0}^{\infty} a_nz^n\in S(g)$, then
\be \label{sub}
T_f(r):=\sum_{n=1}^\infty |a_n|r^n+  \left(\frac{1}{2-\lambda }+\frac{r}{1-r}\right)\|f_0\|_r\leq \lambda
\ee
for $|z|=r<r_g$, and $\lambda =\dist (g(0),\partial g(\ID)) {\color{red}\leq } 1$. The largest such $r_g$, $f\in S(g)$, is called the Bohr radius
in the refined formulation (as described above).

From our earlier two results, we have obtained that Bohr phenomenon in refined formulation exists
for the class of bounded analytic functions and also for the case of analytic functions
with real part less than $1$ in the unit disk. Hence the distance form allows us to extend Bohr's theorem in refined formulation
to a variety of distances. We have the following result which extends Theorem \ref{PVW1-cor1c} in a natural way.
Note that $f(0=g(0)$  and, $f\prec g$ if and only if
$$\frac{f(z)-f(0)}{g'(0)}\prec \frac{g(z)-g(0)}{g'(0)}=z+\frac{1}{g'(0)}\sum_{n=2}^{\infty} \frac{g^{(n)}(0)}{n!}z^n, \quad z\in\ID.
$$
and thus, if needed, it might be convenient to work with normalized superordinate function.

\bthm\label{subtheo1}
Let $f(z)=\sum_{n=0}^{\infty} a_nz^n$ and $g$ be analytic in $\ID$ such that $g$ is univalent and convex in $\ID$.
Assume that $f\in S(g)$ and  $\lambda =\dist (g(0),\partial g(\ID)) {\color{red}\leq } 1$.
Then \eqref{sub}
holds for all  $r\leq r_{*}$, where $r_{*}\approx 0.24683$ as in Theorem  \ref{PVW1-cor1c}.
Moreover, for any  $\lambda \in (0,1)$ there exists a uniquely defined 
$r_0\in \left(r_{*},\frac{1}{3}\right)$ such that $T_f(r)  \le \lambda$
for $r\in [0,r_0]$. The radius $r_0$ is as in Theorem \ref{PVW1-cor1c} given by \eqref{PVW1-eq2-e}. 
\ethm
\bpf
Let $f\prec g$, where  $g(z)=\sum_{n=0}^{\infty} b_nz^n$ is a univalent mapping of $\ID$ onto a convex domain
$\Omega =g(\ID)$.  Then it is well known from the growth estimate for convex functions
and Rogosinksi's coefficient estimate
that (see \cite{DeB1, Rogo-43})
$$\frac{1}{2}|g'(0)|\leq \lambda  \leq |g'(0)|, ~\mbox{ and }~ |b_n| \leq  |g'(0)|~\mbox{ for $n\ge 1$},
$$
where $\lambda =\dist (g(0),\partial \Omega )$. It follows then that $|b_n| \leq 2\lambda$ for $n\geq 1$.
Because $f\prec g$, it follows that $\|f_0\|_r\leq \|g_0\|_r$ for each $0\leq r<1$ and
$$\sum_{n=1}^\infty |a_n|r^n\leq \sum_{n=1}^\infty |b_n|r^n
~\mbox{ for  }~ r \leq \frac{1}{3}.
$$
Combining these two inequalities, we see that the desired conclusion follows if we can show the conclusion for
$T_g(r)$, i.e.,
$$T_g(r)=\sum_{n=1}^\infty |b_n|r^n+  \left(\frac{1}{2-\lambda }+\frac{r}{1-r}\right)\|g_0\|_r\leq \lambda.
$$
Finally, because $|b_n| \leq 2\lambda$ for $n\geq 1$, we have
\beqq
T_g(r)&\leq &2\lambda\sum_{n=1}^\infty r^n+  \left(\frac{1}{2-\lambda }+\frac{r}{1-r}\right)4\lambda^2 \sum_{n=1}^\infty r^{2n}\\
&=& \lambda -\lambda\left[\dfrac{\Phi(\lambda,r)}{(2-\lambda)(1-r)(1-r^2)}\right],
\eeqq
where $\Phi(\lambda,r)$ is as in the proof of Theorem \ref{PVW1-cor1c}. Thus,  $T_g(r) \leq \lambda$ holds
whenever $\Phi(\lambda,r) \geq 0$.
Remaining part of the proof follows from the argument in Theorem \ref{PVW1-cor1c}. The sharpness follows from a
suitable half-plane mapping.
\epf

The idea of this section and Theorem \ref{subtheo1} can be applied to many other situations.
Another instance of this is when $g$ is just univalent in $\ID$ (compare with \cite{Abu} where it is shown
that the sharp radius without the consideration of second term in the expression $T_f(r)$ in \eqref{sub} turns
out to be $3-2\sqrt{2}\approx 0.17157$).

\bthm\label{subtheo2}
Let $g$ be an analytic and univalent function in $\ID$, $f\in S(g)$ and $f(z)=\sum_{n=0}^{\infty} a_nz^n$.
Then the inequality
$$
\sum_{n=1}^\infty |a_n|r^n+  \left(\frac{1}{2-\lambda }+\frac{r}{1-r}\right)\|f_0\|_r\leq \lambda
$$
holds for $|z|=r<r_g $, where $\lambda =\dist (g(0),\partial g(\ID))<1$ and $r_g \approx 0.128445 $ is the unique root of the equation
$$(1-6r+r^2)(1-r)^2(1+r)^3-16r^2(1+r^2)=0
$$
in the interval $(0,1)$. The sharpness of $r_g $ is shown by the Koebe function $f(z)=z/(1-z)^2.$
\ethm
\bpf
Let $f\prec g$, where  $g(z)=\sum_{n=0}^{\infty} b_nz^n$ is a univalent mapping of $\ID$ onto a simply connected domain
$\Omega =g(\ID)$.
Then it is well known from the Koebe estimate and Rogosinksi's coefficient estimate for univalent functions
that (see \cite{DeB1, Rogo-43})
\be\label{eq1-subtheo}
\frac{1}{4}|g'(0)|\leq \lambda  \leq |g'(0)|, ~\mbox{ and }~ |b_n| \leq n |g'(0)| ~\mbox{ for $n\ge 1$},
\ee
where $\lambda =\dist (g(0),\partial \Omega )$.   Also, the first inequality above gives $|b_n| \leq 4n\lambda$ for $n\geq 1$. As in
the proof of Theorem \ref{subtheo1},
we easily have
\beqq
T_g(r)&\leq &4\lambda\sum_{n=1}^\infty nr^n+  \left(\frac{1}{2-\lambda }+\frac{r}{1-r}\right)16\lambda^2 \sum_{n=1}^\infty n^2r^{2n}\\
&=& \lambda-\lambda \left [\frac{(1-r)^2-4r}{(1-r)^2} -\left(\frac{1}{2-\lambda}
+\frac{r}{1-r}\right)\frac{16\lambda r^2(1+r^2)}{(1-r^2)^3}\right ]\\
&=& \lambda -\lambda\left[\dfrac{\Psi(\lambda,r)}{(2-\lambda)(1-r)(1-r^2)^3}\right],
\eeqq
where the equality in the above inequality  is attained when $g(z)$ equals the Koebe function $z/(1-z)^2$, and
\beqq
\Psi(\lambda,r)&=& (1-6r+r^2)(2-\lambda)(1-r)^2(1+r)^3-(1+r(1-\lambda))16\lambda r^2(1+r^2)\\
&=&  16\lambda^2 r^3(1+r^2)- \lambda [(1-6r+r^2)(1-r)^2(1+r)^3 +16r^2(1+r)(1+r^2)]\\
&& \hspace{2cm} +2(1-6r+r^2)(1-r)^2(1+r)^3.
\eeqq
We claim that $\Psi(\lambda,r) \geq 0$ for every $r \leq r_g$ and for $\lambda \in (0,1]$. Clearly,
$$ \frac{\partial^2 \Psi(\lambda,r) }{\partial \lambda^2} \geq 0 ~\mbox{ for every $\lambda \in (0,1]$}
$$
which implies that
\beqq
\frac{\partial \Psi(\lambda,r)}{\partial \lambda} &\leq &  \frac{\partial \Psi}{\partial \lambda}(1,r)\\
& = & -(1-r)[16r^2(1+r^2)+ (1-6r+r^2)(1-r)(1+r)^3] \\
&=& -(1-r)[1-4r+5r^2+27r^4+4r^5-r^6]\\
&=& -(1-r)[r^5(1-r)+r^2+27r^4+3r^5+(2r-1)^2]
\eeqq
from which we obtain that  $\Psi$ is an decreasing function of $\lambda$ on $(0,1)$ so that
$$ \Psi(\lambda, r) \geq \Psi(1,r)=(1-6r+r^2)(1-r)^2(1+r)^3-16r^2(1+r^2)
$$
which is greater than or equal to $0$ for all $r \leq r_g$, where $r_g$ is  as in the statement.
The sharpness of $r_g $ can be easily shown by the Koebe function $f(z)=z/(1-z)^2.$
\epf

\subsection*{Acknowledgments}
The work of the first author is supported by Mathematical Research Impact Centric Support of DST, India  (MTR/2017/000367).


\begin{thebibliography}{99}


\bibitem{Abu} Y. Abu-Muhanna,
Bohr's phenomenon in subordination and bounded harmonic classes,
\emph{Complex Var. Elliptic Equ.} \textbf{55}(11) (2010),  1071--1078.

\bibitem{AAPon1}  Y. Abu-Muhanna, R. M. Ali  and S. Ponnusamy,
On the Bohr inequality,
In ``Progress in Approximation Theory and Applicable Complex Analysis'' (Edited by N.K. Govil et al. ),
Springer Optimization and Its Applications \textbf{117} (2016), 265--295.

%
%
%


\bibitem{AliBarSoly}  R. M. Ali, R. W. Barnard and A. Yu. Solynin,
A note on the Bohr's phenomenon for power series,
\emph{J. Math. Anal. Appl.}  \textbf{449}(1) (2017), 154--167.

\bibitem{AlkKayPon1} S.A. Alkhaleefah, I R. Kayumov, and S. Ponnusamy,
 On the Bohr inequality with a fixed zero coefficient,
\emph{Proc. Amer. Math. Soc.} \textbf{147}(12) (2019), 5263--5274.

\bibitem{BDK5} C. B\'en\'eteau, A. Dahlner and D. Khavinson,
Remarks on the Bohr phenomenon,
\emph{Comput. Methods Funct. Theory} {\bf 4}(1) (2004),  1--19.

\bibitem{BhowDas-98} B.~Bhowmik and N.~Das,
Bohr phenomenon for subordinating families of certain univalent functions,
\emph{J. Math. Anal. Appl.} {\bf 462}(2) (2018), 1087--1098.


\bibitem{BoasKhavin-97-4} H. P. Boas and  D. Khavinson,
Bohr's power series theorem in several variables,
\emph{Proc. Amer. Math. Soc.} \textbf{125}(10) (1997),  2975--2979.

\bibitem{Bohr-14} H. Bohr,
A theorem concerning power series,
\emph{Proc. London Math. Soc.} \textbf{13}(2) (1914), 1--5.

\bibitem{Bom-62}
 E. Bombieri, Sopra un teorema di H. Bohr e G. Ricci sulle funzioni maggioranti delle serie
di potenze,
\emph{Boll. Un. Mat. Ital.} \textbf{17} (3)(1962), 276--282.

\bibitem{BomBor-04} E. Bombieri and J. Bourgain, A remark on Bohr's inequality,
\emph{Int. Math. Res. Not.} \textbf{80} (2004), 4307--4330.

\bibitem{DeB1} L. de Branges,
A proof of the Bieberbach conjecture.
\emph{Acta  Math.} \textbf{154} (1985), 137--152


\bibitem{DeGarM04} A. Defant, S.~R.~Garcia and M. Maestre, Asymptotic estimates for the first and second
Bohr radii of Reinhardt domains,
\emph{J. Approx. Theory} \textbf{128} (2004), 53--68.

\bibitem{DePre06} A. Defant and C. Prengel, Christopher Harald Bohr meets Stefan Banach. Methods in Banach space
theory, 317--339, London Math. Soc. Lecture Note Ser. 337, Cambridge Univ. Press, Cambridge, 2006.


\bibitem{DjaRaman-2000} P.~B.~Djakov and M.~S.~Ramanujan,
A remark on Bohr's theorems and its generalizations,
\emph{J. Analysis} \textbf{8} (2000), 65--77.

\bibitem{DurenUniv-83-8}  P. L. Duren,
Univalent Functions. Springer, New York (1983)

\bibitem{EvPoRa-2017}
S. Evdoridis, S. Ponnusamy and  A. Rasila,
Improved Bohr's inequality for locally univalent harmonic mappings,
\emph{Indag. Math. (N.S.)}, \textbf{30} (2019), 201--213.

\bibitem{GarMasRoss-2018}
S.~R.~Garcia, J.~ Mashreghi and W.~T.~Ross,
\emph{Finite Blaschke products and their connections}, Springer, Cham, 2018.


\bibitem{KayPon1} I. R. Kayumov and S. Ponnusamy,
Bohr inequality for odd analytic functions,
\emph{Comput. Methods Funct. Theory} \textbf{17} (2017), 679--688.

\bibitem{KayPon3} I. R. Kayumov and S. Ponnusamy, Improved version of Bohr's inequality,
\emph{C. R. Math. Acad. Sci. Paris}
\textbf{356}(3) (2018),  272--277

\bibitem{KayPon2} I. R. Kayumov and S. Ponnusamy,
Bohr's inequalities for the analytic functions with lacunary series and harmonic functions,
\emph{J. Math. Anal. and Appl.,}  \textbf{465} (2018), 857--871.



\bibitem{KayPonShak1} I. R. Kayumov, S. Ponnusamy and N. Shakirov,
Bohr radius for locally univalent harmonic mappings,
\emph{Math. Nachr.}  \textbf{291} (2018), 1757--1768.

\bibitem{LiuShangXu-89}
M.~S. Liu, Y.~M. Shang, and J.~F. Xu, Bohr-type inequalities of analytic functions,
\emph{J. Inequal. Appl.} (2018), Paper No. 345, 13 pp

\bibitem{MacGre-67} T. H. MacGregor,  Majorization by univalent functions,
\emph{Duke Math. J.} \textbf{34} (1967), 95--102
%
\bibitem{PaulPopeSingh-02-10} V. I. Paulsen, G. Popescu and D. Singh,
On Bohr's inequality,
\emph{Proc. London Math. Soc. } {\bf 85}(2) (2002), 493--512.

%


\bibitem{PVW1-preprint} S. Ponnusamy, R. Vijayakumar and K.-J. Wirths, {\color{red}
New Inequalities for the Coefficients of unimodular bounded Functions,
Results Math \textbf{75}, 107 (2020). https://doi.org/10.1007/s00025-020-01240-1}



\bibitem{Robert-70} M. S. Robertson, Quasi-subordination and coefficient conjectures,
 \emph{Bull. Amer. Math.Soc.}, \textbf{76} (1970), 1--9.


\bibitem{Rogo-43}  W. Rogosinski, On the coefficients of subordinate functions,
\emph{Proc. London Math. Soc.} \textbf{48}(2) (1943), 48--82.

\bibitem{Sidon-27-15} S. Sidon, \"{U}ber einen Satz von Herrn Bohr,
\emph{Math. Z.} \textbf{26}(1) (1927), 731--732.

\bibitem{Tomic-62-16} M. Tomi\'c, Sur un th\'eor\`eme de H. Bohr,
\emph{Math. Scand.} \textbf{11} (1962), 103--106.

\end{thebibliography}
\end{document}